\documentclass[11pt, a4paper, oneside]{amsart}

\usepackage[utf8]{inputenc}

\usepackage{fullpage}
\usepackage[british]{babel}
\usepackage{amsmath, amsthm, amssymb, mathtools, mathrsfs}
\usepackage{enumerate}
\usepackage{color}
\usepackage{hyperref}

\usepackage{graphicx}
\usepackage{wrapfig}

\newtheorem{theoremcounter}{Theorem Counter}[section]

\theoremstyle{definition}

\newtheorem{remark}[theoremcounter]{Remark}

\theoremstyle{plain}

\newtheorem{theorem}[theoremcounter]{Theorem}

\newcommand{\Z}{\mathbb{Z}}

\newcommand{\calP}{\mathcal{P}}

\newcommand{\scB}{\mathscr{B}}

\DeclareMathOperator{\Li}{Li}

\newcommand{\mat}[1]{\begin{matrix}#1\end{matrix}}

\newcommand{\bmat}[1]{\begin{bmatrix}#1\end{bmatrix}}
\newcommand{\smat}[1]{\bigl[\begin{smallmatrix}#1\end{smallmatrix}\bigr]}

\newcommand{\thmref}[2]{\hyperref[#2]{#1 \ref*{#2}}}
\renewcommand{\eqref}[1]{\hyperref[#1]{(\ref*{#1})}}

\renewcommand{\emph}[1]{\textbf{#1}}

\begin{document}

\title{Symmetrized poly-Bernoulli numbers and combinatorics}
\author{Toshiki Matsusaka} 
\date{\today}

\address{Faculty of Mathematics, Kyushu University, Motooka 744, Nishi-ku Fukuoka 819-0395, Japan,}
\email{toshikimatsusaka@gmail.com}


\thanks{This work was supported by JSPS KAKENHI Grant Numbers 18J20590.}

\maketitle

\begin{abstract}
	Poly-Bernoulli numbers are one of generalizations of the classical Bernoulli numbers. Since a negative index poly-Bernoulli number is an integer, it is an interesting problem to study this number from combinatorial viewpoint. In this short article, we give a new combinatorial relation between symmetrized poly-Bernoulli numbers and the Dumont-Foata polynomials.
\end{abstract}



\section{Introduction}


A poly-Bernoulli polynomial $B_m^{(\ell)} (x)$ of index $\ell \in \Z$ is defined by the generating series
\[
	e^{-xt} \frac{\Li_\ell (1 - e^{-t})}{1 - e^{-t}} = \sum_{m=0}^\infty B_m^{(\ell)} (x) \frac{t^m}{m!},
\]
where $\Li_\ell (z)$ is the polylogarithm function given by
\[
	\Li_\ell (z) = \sum_{m=1}^\infty \frac{z^m}{m^\ell} \quad (|z| < 1).
\]
The polynomial $B_m^{(1)}(x)$ coincides with the classical Bernoulli polynomial $B_m(1-x) = (-1)^m B_m(x)$ since $\Li_1(z) = -\log(1-z)$ holds (see \cite[(4.2)]{AIK14}). Following Kaneko \cite{Kan97}, the special value $B_m^{(\ell)} := B_m^{(\ell)}(0)$ is particularly called a poly-Bernoulli number of index $\ell$. \\

The aim of this study is to give a combinatorial perspective to the special values of $B_m^{(\ell)}(x)$ at integers $k \in \Z$. We now assume that the index $\ell \leq 0$, so that the values $B_m^{(\ell)}(k)$ are always integers. One of the first such investigations is the coincidence of two numbers, the poly-Bernoulli number $B_m^{(\ell)}(0)$ and the number of $01$ lonesum matrices of size $m \times |\ell|$ (see \cite{Bre08}). A recent work by B\'{e}nyi-Hajnal \cite{BH17} established more combinatorial relations in this direction.\\ 

In this article, we take a step in another direction along with Kaneko-Sakurai-Tsumura \cite{KST18}. To describe this more precisely, let $G_n$ be the Genocchi number defined by $G_n = 2(2^{n+2} - 1)|B_{n+2}|$, where $B_m = B_m^{(1)}(1)$ is the classical Bernoulli number. Then they showed in \cite[Theorem 4.2]{KST18} that 
\begin{align}\label{k=1}
	\sum_{\ell =0}^n (-1)^\ell B_{n-\ell}^{(-\ell -1)} (1) = (-1)^{n/2} G_n
\end{align}
for any $n \geq 0$. As they mentioned, this equation is an analogue of Arakawa-Kaneko's result \cite{AK99}
\begin{align}\label{k=0}
	\sum_{\ell =0}^n (-1)^\ell B_{n-\ell}^{(-\ell)} (0) = \left\{\begin{array}{ll}
		1 &\text{if } n = 0\\
		0 &\text{if } n > 0.
	\end{array} \right.
\end{align}
In addition, Sakurai asked in her master thesis whether we can generalize these equations for any positive integers $x = k$, and give some combinatorial meaning to them. Our main result provides an answer to these two questions in terms of the Dumont-Foata polynomial as follows.

\begin{theorem}\label{main}
	Let $G_n(x,y,z)$ be the $n$-th Dumont-Foata polynomial defined in $(\ref{DF-poly})$, and $\scB_m^{(\ell)}(k)$ the symmetrized poly-Bernoulli number defined in $(\ref{SPB})$. Then we have
	\[
		\sum_{\ell=0}^n (-1)^\ell \scB_{n-\ell}^{(-\ell)} (k) = k! \cdot (-1)^{n/2} G_n(1,1,k)
	\]
	for any non-negative integers $n, k \geq 0$. In particular, both sides equal zero for odd $n$.
\end{theorem}

This theorem recovers the equations (\ref{k=1}) and (\ref{k=0}) since $\scB_m^{(-\ell)} (0) = B_m^{(-\ell)} (0), \scB_m^{(-\ell)} (1) = B_m^{(-\ell-1)}(1)$, and $G_n(1,1,1) = G_n$ hold as we see later. 

\begin{remark}
Recent work of B\'{e}nyi-Hajnal \cite[Section 6]{BH17} pointed out that the sequence of diagonal sums
\[
	\sum_{\ell = 0}^n B_{n-\ell}^{(-\ell -1)} (1) \quad (n \geq 0)
\]
appears in OEIS \cite{OEIS} A136127. It might be an interesting problem to extend our theorem in this direction.
\end{remark}


\section{Definitions}



\subsection{Dumont-Foata polynomials}


We review the work of Dumont-Foata \cite{DF76} here. For a positive even integer $n \in 2\Z$, we consider the surjective map $p: \{1,2,\dots, n\} \to \{2,4, \dots, n\}$ with $p(x) \geq x$ for each $x \in \{1,2, \dots, n\}$. This map is called (surjective) pistol of size $n$, and corresponds to the following diagram. Here we draw one example for $n = 6$ given by $p(1) = 2, p(2) = 4, p(3) = 6, p(4) = 4, p(5) = 6, p(6) = 6$.

\begin{figure}[h]
	\includegraphics[width=50mm]{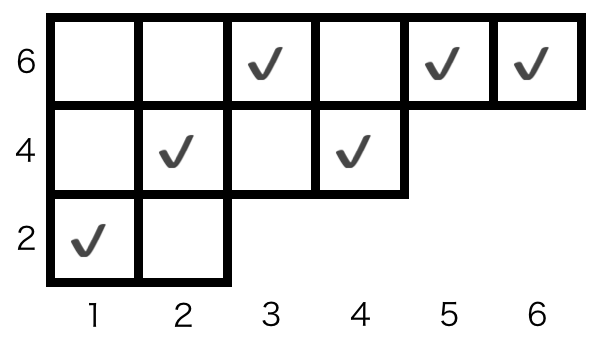} 
\end{figure}

Let $\calP_n$ be the set of all pistols of size $n$. For each pistol $p \in \calP_n$, we define three quantities called bulging, fixed, and maximal points. First, the number $x \in \{1, 2, \dots, n\}$ is a bulging point of $p \in \calP_n$ if $p(y) < p(x)$ for any $0< y < x$. Then we denote by $b(p)$ the number of bluging points of $p$. In the diagram, $b(p)$ corresponds to the number of steps of the minimal stair covering all check marks. For the above example $p$, the points $x = 1, 2, 3$ are bulging points, so that $b(p) = 3$.

\begin{figure}[h]
	\includegraphics[width=100mm]{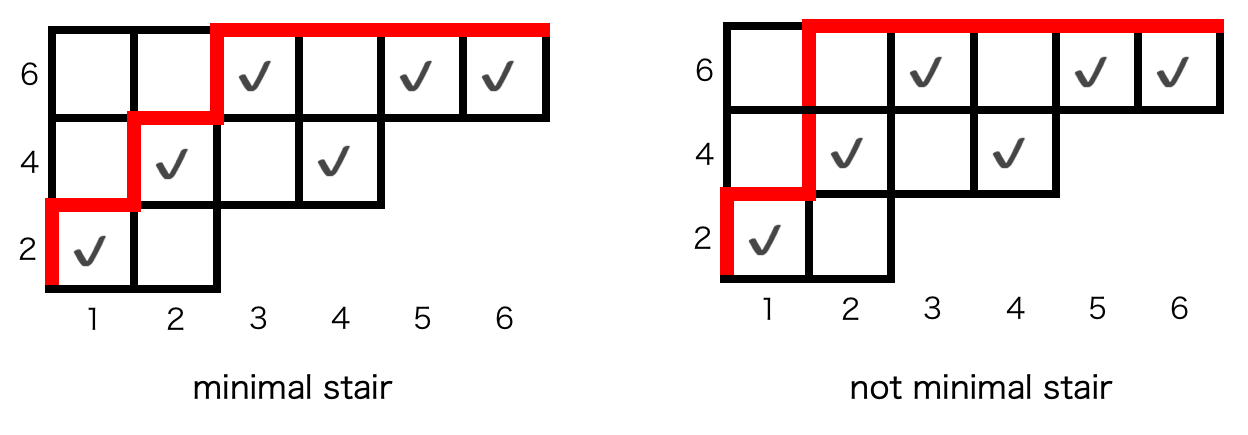} 
\end{figure}

Next, the point $x \in \{1, 2, \dots, n\}$ is called a fixed point of $p \in \calP_n$ if $p(x) = x$. Finally, the point $x \in \{1, 2, \dots, n-1\}$ is a maximal point of $p \in \calP_n$ if $p(x) = n$. Then we denote by $f(p)$ and $m(p)$ the numbers of fixed points and maximal points of $p \in \calP_n$, respectively. For the above example, we have $f(p) = 2$ and $m(p) = 2$. Under these notations, Dumont-Foata established the following interesting theorem.

\begin{theorem}\cite[Th\'{e}or\`{e}me 1a, 2]{DF76}
	Let $n \in 2\Z_{>0}$ be a positive even integer. The polynomial defined by
	\begin{align}\label{DF-poly}
		G_n(x,y,z) := \sum_{p \in \calP_n} x^{b(p)} y^{f(p)} z^{m(p)}
	\end{align}
	is a symmetric polynomial in three variables, and gives the Genocchi number $G_n$ as
	\[
		G_n(1,1,1) = G_n.
	\]
\end{theorem}

In addition, we put $G_0(x,y,z) = 1$ and $G_n(x,y,z) = 0$ for a positive odd integer $n \in \Z$. The polynomial $G_n(x,y,z)$ is called the $n$-th Dumont-Foata polynomial. Furthermore, they showed that the polynomial for $n > 0$ has the form $G_n(x,y,z) = xyz F_n(z,y,z)$, and the polynomial $F_n(x,y,z)$ satisfies the recurrence relation
\[
	F_n(x,y,z) = (x+z)(y+z) F_{n-2}(x,y,z+1) - z^2 F_{n-2}(x,y,z)
\]
with initial values $F_1(x,y,z) = 0$ and $F_2(x,y,z) = 1$. This implies that the polynomial $G_n(z) := G_n(1,1,z)$ called the Gandhi polynomial satisfies
\begin{align}\label{RF}
	G_{n+2}(z) = z(z+1) G_n(z+1) - z^2 G_n(z)
\end{align}
with $G_0(z) = 1, G_1(z) = 0$. \\

For instance, there exist three pistols of size $4$. Each pistol has $(b(p), f(p), m(p)) = (2,2,1), (2,1,2)$ and $(1,2,2)$, so that the Dumont-Foata polynomial is given by
\[
	G_4(x,y,z) = x^2 y^2 z + x^2 y z^2 + x y^2 z^2 = xyz(xy + yz + zx).
\]
In fact, $G_4(1,1,1) = 3$ coincides the $4$-th Genocchi number given by $G_4 = 2(2^6 -1) |B_6| = 126 \times 1/42 = 3$.

\begin{figure}[h]
	\includegraphics[width=110mm]{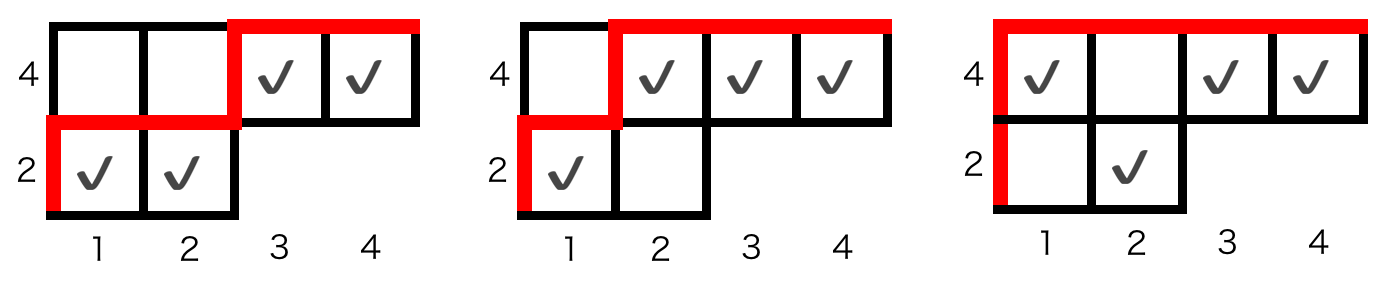} 
\end{figure}


\subsection{Symmetrized poly-Bernoulli numbers}


The following are the tables of poly-Bernoulli numbers $\{B_m^{(-\ell)}(0)\}$ and $\{B_m^{(-\ell)}(1)\}$ with $m, \ell \geq 0$.

\begin{table}[htbp]
	\caption{$B_m^{(-\ell)} (0)$: left, $B_m^{(-\ell)} (1)$: right}
	\begin{tabular}{|c||c|c|c|c|c|c} \hline
		$\ell \backslash m$ & 0 & 1 & 2 & 3 &4\\ \hline\hline
		0 & 1 & 1 & 1 & 1 & 1\\ \hline
		1 & 1 & 2 & 4 & 8 & 16\\ \hline
		2 & 1 & 4 & 14 & 46 & 146\\ \hline
		3 & 1 & 8 & 46 & 230 & 1066\\ \hline
		4 & 1 & 16 & 146 & 1066 & 6902\\ \hline
	\end{tabular}
	\hspace{20mm}
	\begin{tabular}{|c||c|c|c|c|c|c} \hline
		$\ell \backslash m$ & 0 & 1 & 2 & 3 &4\\ \hline\hline
		0 & 1 & 0 & 0 & 0 & 0\\ \hline
		1 & 1 & 1 & 1 & 1 & 1\\ \hline
		2 & 1 & 3 & 7 & 15 & 31\\ \hline
		3 & 1 & 7 & 31 & 115 & 391\\ \hline
		4 & 1 & 15 & 115 & 675 & 3451\\ \hline
	\end{tabular}
\end{table}

\noindent We can find the symmetric properties of these numbers at a glance. On the other hand, for $k \geq 2$ it seems unlikely that such a simple symmetric property can be given. In order to reproduce the symmetric properties for any $k \geq 2$, Kaneko-Sakurai-Tsumura \cite{KST18} considered combinations of $B_m^{(-\ell)} (k)$. To make this precise, let $m, \ell, k \geq 0$ be non-negative integers. We now define the symmetrized poly-Bernoulli number $\scB_m^{(-\ell)} (k)$ by
\begin{align}\label{SPB}
	\scB_m^{(-\ell)} (k) = \sum_{j=0}^k \bmat{k \\ j} B_m^{(-\ell - j)} (k),
\end{align}
where $\smat{k \\ j}$ is the Stirling number of the first kind defined in \cite[Definition 2.5]{AIK14}. Then this number satisfies $\scB_m^{(-\ell)} (k) = \scB_\ell^{(-m)} (k)$ for any $m,\ell, k \geq 0$, and
\[
	\scB_m^{(-\ell)} (0) = B_m^{(-\ell)} (0), \quad \scB_m^{(-\ell)} (1) = B_m^{(-\ell-1)} (1).
\]

\begin{table}[htbp]
	\caption{$B_m^{(-\ell)} (2)$: left, $\scB_m^{(-\ell)} (2)$: right}
	\begin{tabular}{|c||c|c|c|c|c|c} \hline
		$\ell \backslash m$ & 0 & 1 & 2 & 3 &4\\ \hline\hline
		0 & 1 & $-1$ & 1 & $-1$ & 1\\ \hline
		1 & 1 & 0 & 0 & 0 & 0\\ \hline
		2 & 1 & 2 & 2 & 2 & 2\\ \hline
		3 & 1 & 6 & 18 & 42 & 90\\ \hline
		4 & 1 & 14 & 86 & 374 & 1382\\ \hline
	\end{tabular}
	\hspace{20mm}
	\begin{tabular}{|c||c|c|c|c|c|c} \hline
		$\ell \backslash m$ & 0 & 1 & 2 & 3 &4\\ \hline\hline
		0 & 2 & 2 & 2 & 2 & 2\\ \hline
		1 & 2 & 8 & 20 & 44 & 92\\ \hline
		2 & 2 & 20 & 104 & 416 & 1472\\ \hline
		3 & 2 & 44 & 416 & 2744 & 15032\\ \hline
		4 & 2 & 92 & 1472 & 15032 & 120632\\ \hline
	\end{tabular}
\end{table}

Moreover they showed the following explicit formula for $\scB_m^{(-\ell)} (k)$.
\begin{align}\label{EF}
	\scB_m^{(-\ell)} (k) = \sum_{j=0}^{\min(m, \ell)} j! (k+j)! \left\{ \mat{m + 1 \\ j+1} \right\} \left\{ \mat{\ell + 1 \\ j+1} \right\},
\end{align}
where $\left\{\begin{smallmatrix}k \\ j\end{smallmatrix} \right\}$ is the Stirling number of the second kind defined in \cite[Definition 2.2]{AIK14}. By using this formula, we prove our main theorem in the next section.


\section{Proof}


To prove Theorem \ref{main}, it suffices to show that the function
\[
	\tilde{G}_n(k) := \frac{(-1)^{n/2}}{k!} \sum_{\ell=0}^n (-1)^\ell \scB_{n-\ell}^{(-\ell)} (k)
\]
satisfies the recurrence relation (\ref{RF}) for any integer $k \geq 0$. First we can easily see that $\tilde{G}_0(k) = 1$ and $\tilde{G}_1(k) = 0$, which are initial cases. Moreover for any odd integer $n$,  $\tilde{G}_n(k) = 0$ follows from the symmetric property of $\scB_m^{(-\ell)}(k)$. For an even integer $n \geq 2$, by the formula (\ref{EF}) we have

\begin{align}
	(-1)^{n/2} &k! \left[k(k+1) \tilde{G}_n (k+1) - k^2 \tilde{G}_n (k) - \tilde{G}_{n+2} (k) \right] \nonumber\\
		&= k \sum_{j=0}^{n/2} j! (k+j+1)!  \sum_{\ell = j}^{n - j} (-1)^\ell \left\{ \mat{n-\ell + 1 \\ j+1} \right\} \left\{ \mat{\ell + 1 \\ j+1} \right\} \label{1}\\
		&\qquad -k^2 \sum_{j = 0}^{n/2} j! (k+j)! \sum_{\ell = j}^{n-j} (-1)^\ell \left\{ \mat{n-\ell + 1 \\ j+1} \right\} \left\{ \mat{\ell + 1 \\ j+1} \right\} \label{2}\\
		&\qquad + \sum_{j=0}^{n/2+1} j! (k+j)! \sum_{\ell = j}^{n+2-j} (-1)^\ell \left\{ \mat{n-\ell +3 \\ j+1} \right\} \left\{ \mat{\ell + 1 \\ j+1} \right\}. \label{3}
\end{align}
Since $\left\{\begin{smallmatrix}k \\ 1\end{smallmatrix} \right\} = 1$ holds for any $k \geq 1$, the third line (\ref{3}) is divided according as $j = 0$ or not, which equals
\[
	k! + \sum_{j=0}^{n/2} (j+1)! (k+j+1)! \sum_{\ell = j +1}^{n+1-j} (-1)^\ell \left\{ \mat{n +3-\ell \\ j+2} \right\} \left\{ \mat{\ell+1 \\ j +2} \right\}.
\]
Let
\begin{align}\label{a}
	a_{n,j} := \sum_{\ell = j}^{n-j} (-1)^\ell \left\{ \mat{n-\ell + 1 \\ j+1} \right\} \left\{ \mat{\ell + 1 \\ j+1} \right\},
\end{align}
then the total of (\ref{1}), (\ref{2}), and (\ref{3}) equals
\begin{align}\label{goal}
	k! + \sum_{j=0}^{n/2} j!(k+j)! \bigg[ k(j+1) a_{n,j} + (j+1)(k+j+1) a_{n+2, j+1} \bigg].
\end{align}
Once the sum is $0$, the proof completes. By using the generating function given in \cite[Proposition 2.6 (8)]{AIK14}, we have
\begin{align*}
	\frac{t^j}{(1-t)(1-2t) \cdots (1-(j+1)t)} = \sum_{\ell \geq j} \left\{ \mat{\ell+1 \\ j+1} \right\} t^\ell = \sum_{\ell \leq n - j} \left\{ \mat{n-\ell +1 \\ j+1} \right\} t^{n-\ell}
\end{align*}
for any non-negative integers $n, j \in \Z_{\geq 0}$. Multiplying these two expressions, we obtain
\[
	\frac{s^j t^j}{(1-s)(1-t) \cdots (1-(j+1)s)(1-(j+1)t)} = \sum_{\ell \geq j} \sum_{k \leq n - j} \left\{ \mat{\ell +1 \\ j +1} \right\} \left\{ \mat{n-k+1 \\ j+1} \right\} s^\ell t^{n-k}.
\]
By specializing at $s = -x, t = x$,
\begin{align}\label{gen}
	\frac{(-1)^j x^{2j}}{(1-x^2) \cdots (1-(j+1)^2 x^2)} = \sum_{\ell \geq j} \sum_{k \leq n -j} (-1)^\ell \left\{ \mat{\ell +1 \\ j +1} \right\} \left\{ \mat{n-k+1 \\ j+1} \right\} x^{n+\ell -k}.
\end{align}
Thus we see that the number $a_{n,j}$ defined in (\ref{a}) appears as the $n$-th coefficient of (\ref{gen}). By the expression of the left-hand side of (\ref{gen}), we easily see that $a_{n,j} = 0$ when $n$ is an odd integer or $2j > n$. Further we get the initial values $a_{2j, j} = (-1)^j, a_{n,0} = 1$ for even $n$, and the recurrence relation
\[
	a_{n+2,j} = (j+1)^2 a_{n,j} - a_{n,j-1}.
\]
Applying this to the term in (\ref{goal}), we get
\begin{align*}
	k! + \sum_{j=0}^{n/2} &j!(k+j)! \bigg[ k(j+1) a_{n,j} + (j+1)(k+j+1) ((j+2)^2 a_{n, j+1} - a_{n, j}) \bigg]\\
		&= k! + \sum_{j=0}^{n/2} (j+2)(j+2)!(k+j+1)! a_{n,j+1} - \sum_{j=0}^{n/2} (j+1) (j+1)! (k+j)! a_{n,j}.
\end{align*}
Since $a_{n, n/2+1} = 0$ and $a_{n,0} = 1$, this equals $0$, which concludes the proof of Theorem  \ref{main}.


\end{document}